\documentclass[12pt,a4,twoside]{amsart}
\usepackage[english]{babel} 
\usepackage{amssymb}
\usepackage{amsmath}
\usepackage{txfonts}
\usepackage{mathdots}
\usepackage[classicReIm]{kpfonts}
\usepackage[dvips]{graphicx} 

\theoremstyle{plain}
\newtheorem{theorem}{Theorem}

\newtheorem{proposition}{Proposition}

\theoremstyle{definition}

\theoremstyle{plain}
\newtoks\thehProclaim
\newtheorem*{Proclaim}{\the\thehProclaim}

\theoremstyle{remark}
\newtheorem{remark}{\sc  Remark}
\theoremstyle{remark}

\theoremstyle{definition}
\newtoks{\thehRemark}
\newtheorem*{Remark}{\the\thehRemark}

\begin{document}


\title{CONNECTING CYCLES FOR CONCENTRIC CIRCLES}

\author{George Khimshiashvili, Dirk Siersma}

\address{Ilia State University, Tbilisi;  Utrecht University}

\email{gogikhim@yahoo.com ;  d.siersma@uu.nl}

\begin{abstract}We study perimeters of connecting cycles for concentric circles. More precisely, we are interested in characterization of those connecting cycles which are critical points of perimeter considered as a function on the product of given circles. Specifically, we aim at showing that, generically, perimeter is a Morse function on the configuration space, and computing Morse indices of critical configurations. In particular, we prove that the diametrically aligned configurations are critical and their indices can be calculated from an explicitly given tridiagonal matrix. For four concentric circles, we give examples of non-generic collections of radii and describe a pitchfork type bifurcation of stationary connecting cycles.
\end{abstract}



\maketitle


\section{Introduction}

\noindent 

\noindent The problem of geometrical characterization and construction of \textit{minimal connecting cycles} for a system of non-intersecting planar domains or contours arises in the framework of various topics of combinatorial geometry and applied mathematics \cite{bdy}, \cite{pa}. In the present paper we consider this problem in a wider context of extremal problems and Morse theory on configuration spaces. To explain and illustrate our approach we discuss it in some detail for systems of three and four concentric circles. We also give several generalizations and indicate some similar problems suggested by our approach.

\medskip

\noindent We begin with describing the aims and strategy of research in the case of three concentric circles. In particular, we prove that, generically, perimeter defines a non-degenerate (Morse) function on the reduced configuration space, find the number of its critical points and compute their indices. Next, we investigate the same issues for four concentric circles and show that several new interesting phenomena appear in this case, which suggests some general conjectures and lines of further research. In conclusion, we briefly discuss several possible generalizations and related problems.

\medskip

\noindent The main results presented in this paper have been obtained in the framework of a "Research in Pairs" project at the Centre de Rencontres Mathematique (CIRM) in Luminy, France in April of 2018 and were finalized at a ``Research in Groups'' project at the International Centre for Mathematical Sciences (ICMS) in Edinburgh, Scotland in September 2018. The authors wish to acknowledge the hospitality and excellent working conditions at CIRM and ICMS. They thank Gaiane Panina for useful discussions and contributions to the project.

\section{Connecting cycles and  equilibrium configurations}

\noindent We proceed by describing a general setting for our approach and particular cases studied in the sequel. Consider a system $\Gamma$ of $n$ non-intersecting smooth contours (simple closed differentiable curves)  $\gamma_i$ in Euclidean plane   $\mathbb{R}^2$. An ordered collection of $n$ points $p_i \in \gamma_i$  defines an $n$-gon $P$  with the sides  $p_i p_{i+1}$. Such polygons will be called \textit{connecting cycles} for$\Gamma$ or $\Gamma$-circuits. The set of all such polygons  $P(\Gamma)$ is naturally identified with the Cartesian product $\prod  \gamma _i$ of given contours. We endow $P$ with the smooth structure inherited from    $\prod  \gamma _i$  so that it becomes diffeomorphic to $n$-torus \textbf{$\mathbb{T}$}\textit{${}^{n}$}.                               

\medskip 

\noindent Perimeter $L(P)$ of a $\Gamma$-circuit $P$ defines a smooth function on $P(\Gamma)$. A classical problem of combinatorial geometry is to find a point of global minimum of $L$ on $p(\Gamma)$which is called a \textit{minimal connecting cycle} [5]. We extend this problem by considering all critical points of $L$  referred to as \textit{stationary connecting cycles} for $\Gamma$ or\textit{ stationary $\Gamma$-circuits}. It seems rather hopeless to obtain useful description of stationary circuits in general setting, but it turns out that in several special cases our approach leads to interesting constructions and results part of which are presented in the sequel.

\begin{remark}
 A closely related topic is concerned with the so-called minimal connections for a system of non-intersecting compact sets \cite{dt}]. We emphasize that the problem of minimal connections is essentially different from our setting. In particular, a minimal connection need not be a part of a minimal connecting cycle [3]. We do not discuss minimal connections in this paper.
\end{remark}

\noindent In the sequel it is convenient to think of stationary circuits as \textit{equilibrium configurations. }The latter term is motivated by the following informal but useful "mechanical interpretation" of stationary circuits. Imagine that instead of points we have a system of very small rigid circles $\delta_ i$  such that $\delta_i$  is linked with $\gamma_i$ and can freely slide along  $\gamma_i$  with its centre remaining on $\gamma_i$. Suppose further that we also have a sufficiently short elastic closed string in  satisfying Hook's law, which is linked with each $\delta_ i$  ("sufficiently short" in this context means that the string is stretched for any $\Gamma$-circuit, i.e., is shorter than the minimal $\Gamma$-circuit). Consider now the equilibrium configurations of such a mechanical system assuming that positions of contours $\Gamma_i$   remain fixed, i.e., we are interested in those positions of circles  $\delta_ i$  in which they stay in rest. By a general law of mechanics this happens if and only if the elastic energy of the string takes one of its stationary values. Since by Hook's law the elastic energy is a constant multiple of the squared length of the string it follows that the equilibrium configurations are the same as the stationary circuits.

\medskip

\noindent Using this interpretation, one notices that, for an equilibrium configuration, the resultant force at each circle $\delta_i$  should be collinear with the radius connecting this point with the centre of $\delta_i$. As the radii of circles $\delta_i$  tend to zero, it follows that the two adjacent sides of stationary $\Gamma$-circuit are either bisected by the radius to this point or lie on the same straight line. This property resembles the famous \textit{Fermat principle} for light rays \cite{bdy}. From this point of view the case where two sides are continuation of each other can be considered as refraction. We will show that this analogy leads to useful conclusions and refer to such polygons as \textit{Fermat $\Gamma$}-circuits. 

\medskip

\noindent From now on we basically deal with the case where $\Gamma$  is a system  $S= \{ C(r_i)\}$ of $n  \ge 3$ coplanar concentric circles with given radii $r_i$. To avoid consideration of non-isolated critical points, in the sequel we always fix a point on the outer circle and work in the reduced configuration space $M'= M'(S)$,
  diffeomorphic to \textit{(n-1)}-torus $\mathbb{T}^{n-1}$ equipped with a system of natural coordinates given by the polar angles of vertices.$M'$ is in fact the configuration space of two concentric circles and an external point. The original circle symmetry of the problem is now reduced to a reflection in a line. Solutions will occur in pairs unless they are completely aligned.

\noindent 

\noindent In this setting, our aims can be described as:
\begin{itemize}
\item[(1)]  investigating if $L$ is a Morse function on $M'$;
\item[(2)]  finding the maximal possible number of critical points of $L$ on $M'$;
\item[(3)]  calculating the number of critical points for concrete values of radii;
\item[(4)]  calculating the Morse indices of non-degenerate critical points of  $L$;
\item[(5)] investigating bifurcations and degenerate critical points arising in this context. 
\end{itemize}

\noindent In relation with the last two problems it appears instructive to pay special attention to a special type of centrally aligned configurations called \textit{parades }by a way of analogy with the term\textit{ parades of planets }used in celestial mechanics. It turns out that parades are stationary circuits for arbitrary $n$ which follows, in particular, from our first main result. 

\begin{theorem}

\begin{itemize}

\item[(a)]
 Critical points of L on $M'(S)$ are exactly the Fermat S-circuits,
\item[(b)]
For each stationary circuit P, there exists a circle, concentric with the circles of S, such that each side of  P or its straight line continuation is tangent to that circle.
\end{itemize}
\end{theorem}
\noindent 

\noindent The proof uses Lagrange multipliers and elementary geometric considerations. The circle mentioned in the theorem will be called \textit{tangential circle }of circuit  $P$.

\begin{remark}
Simple examples show that stationary circuits may have various
shapes. For our purposes it is sufficient to distinguish four types: convex,
non-convex without self-intersections, self-intersecting, and (partially) aligned.
Parades provide examples of stationary circuits which exhibit combinations of refracted (aligned) and reflected cases. As we will see, there may also exist partially aligned stationary circuits different from parades.
\end{remark}

\noindent Each Fermat circuit has a well-defined sequence of reflections and refractions. Its tangential circle mentioned in Theorem 1 has a common centre with the circles of S and its radius is smaller than any of the radii of given concentric circles. In the case of parades, we accept the convention that the tangential circle consists of the centre only and its radius is $0$. A stationary circuit, which exhibits only reflections is called it a \textit{Snellius circuit}. In the sequel we refer to its tangential circle as the \textit{socle }of stationary circuit. It is also easy to see that, for each socle of radius $\sigma$, there exist exactly two (mirror symmetric) stationary circuits circumscribing it and their perimeters are equal to 
$2 \sum{\sqrt{r_i^2 - \sigma_i^2}}$ . Thus knowing all possible values of radii of socles one can estimate the number of stationary circuits of Snellius type. So finding all possible radii of socles is another natural problem which is also addressed in the sequel.

\medskip 

\noindent We take now a closer look at a special class of stationary circuits mentioned above as parades. More precisely, \textit{parade} is a $S$-circuit $P$ such that the polar angles of all of its vertices are either 0 or $\pi$. Obviously, the number of parades is  $2^n$. Let us once and forever choose a coordinate system with the common centre of given circles at the origin and fix $p_n$ as the intersection point of the outer circle with the positive ray on  $x$-axis. The number of parades with fixed $p_n$ is $2^{n-1}$. Each parade obviously satisfies the Fermat principle so by Theorem 1 it is a stationary $S$-circuit. This can also be seen from the following explicit formula for the gradient of perimeter which will also be used for investigating non-degeneracy of parades. To present the aforementioned formula we use polar coordinates and write  $Op_i = r_i(\cos \alpha_i,\sin\alpha_i)$
 Note that by our assumption $\alpha_n=0 $. 
For each\textit{ j (mod n)},  we have:
$$ l_j = |p_jp_{j+1}| = \sqrt{r_j^2 + r_{j+1}^2 - 2  r_j r_{j+1} \cos ( \alpha_{j+1} - \alpha_j)}.
$$
Then by law of cosines the perimeter of $P$ equals $\sum_{j=1}^n l_j$
 and the formula for gradient can be written as follows:

\noindent 
\begin{proposition}
The gradient of perimeter at $S$-circuit  $P$ as above is given by
$$ \left(
 \frac{r_{j-1} r_j \sin (\alpha_j - \alpha_{j-1})}{l_{j-1}} +  \frac{r_{j+1} r_j \sin (\alpha_j - \alpha_{j+1})}{l_{j}}
\right),  j= 1, \cdots, n.
$$ 
\end{proposition}

Note that both  areas and lengths occur in this formula. Using this formula one can compute the Hessian matrix of $L$ for arbitrary $n$ and verify that it has a tridiagonal form and its determinant gives rise to rational expressions in the radii (for details see n=3 and n=4 below). This determinant can be explicitly evaluated at parades, which gives an expression like the hessian formula given in \cite{sw} and yields our second main result.

\begin{theorem}
For generic values of radii $r_i$, all parades are non-degenerate critical points of $L$ on $M'(S)$.
\end{theorem}

\noindent The condition of genericity is essential. For $n=3$, it is equivalent to the requirement that the radii are pairwise distinct. However this condition is not sufficient for bigger $n$.  Namely, we show below that, for $n \ge 4$ , there exist collections of radii such that some of the parades are degenerate critical points of perimeter. Before passing to these subtle general issues we give more details in the case of three concentric circles (\textit{3CCcase}).

\section{Three Concentric Circles (3CC)}

\noindent For  $n=3$, consider three concentric circles having radii
$r_1 < r_2 < r_3$ . First, we show that, in this case, all parades are non-degenerate critical points of the perimeter, and calculate their Morse indices.

\noindent 

\begin{proposition}
 For $n=3$ and pairwise different radii, all parades are non-degenerate critical points of $ L$ on $M'$. The shortest one is a non-degenerate global minimum and the other three parades are non-degenerate saddles (of index one).
\end{proposition}
\noindent 

\begin{proof} Since all vertices are on the same straight line we will assume that this is the -axis and that the point on the third (outer) circle has positive $x$-coordinate.  Then there are four different parades. Let now $x_i$ be the $x$-coordinate of the \textit{i}-th vertex.  Then we have $x_i= \pm  r_i$ if  $i=1, 2$,  and $x_3 = r_3$. Using Proposition 1 one is able to calculate the second partial derivatives and obtain explicit formulae for hessian at any parade. Namely, the Hessian matrix has the form
$$ \left( \begin{array}{cc}
b_1 + b_2 & - b_2              \\
-b_2            & b_2 + b_3 
\end{array}
\right)
$$
\noindent
where $b_i = \frac{x_{i-1}x_i}{|x_{i-1} - x_i|}, \;  i \; \mbox{modulo} \; 3$.

\noindent Now it is easy to verify that, for pairwise distinct radii, the hessian is non-zero at any of the four parades. With the Hessian matrix at hand it is also easy to compute the Morse index of perimeter at any parade using the classical Sylvester rule (we give more details in the 4CC case). The combination of signs in the Hessian matrix depends on the choice of parade, i.e., on the polar angles of the parade points. For example, for the "shortest parade", all polar angles are equal to $0$ and all elements in the above matrix are positive. The determinant is also positive. So by Sylvester rule the Morse index vanishes and this is a non-degenerate minimum as is also clear by geometric considerations. In the same way it is easy to verify that the other three parades are non-degenerate saddles if radii are pairwise different. 
\end{proof}

\noindent We can now give a rather detailed description of stationary connecting cycles in the 3CC case.

\begin{theorem}
For any three distinct concentric circles, perimeter is a Morse function on M' with exactly six critical points: four parades (the short one is the unique minimum, i.e., a genuine minimal connecting cycle, the other three are. saddles), and two non-degenerate triangles, where L attains the global maximum The values of L at all critical points given in increasing order are:
$$ L_m = 2(r_3-r_1) \; ; \; 2(r_1+r_3) \; ; \; 2(r_2+r_3) \; ; $$
$$ L_M = \sqrt{r_1^2+ r_2^2 + 2 t r_1 r_2 r_3^{-1}} +  \sqrt{r_1^2+ r_3^2 + 2 t r_1 r_3 r_2^{-1}} +  \sqrt{r_2^2+ r_3^2 + 2 t r_2 r_3 r_1^{-1}}.
$$
\noindent where  $t$  is the unique positive real root of cubic equation
$$ 2 r_1 r_2  r_3 t^3 + (r_1^2 r_2^2 + r_2^2 r_3^2 + r_3^2 r_1^2) t^2 - r_1^2 r_2 ^2 r_3^2 = 0. $$
\end{theorem}

To show that there are only six critical points we use Lagrange multipliers and elementary algebraic manipulations. Using Fermat principle it is easy to construct the two stationary triangles and get a cubic equation for the maximum of perimeter. The types of other critical points are already known from Proposition 2.

\begin{remark}
Notice that L is not an exact Morse function, since we have 6 critical points and the Betti numbers of the moduli space add up to 4. Also note that the 
critical values  $2(r_2 + r_3)$ and $L_M$ are double, i.e., there are two critical points on each of those levels (For $L_M$ this comes from the reflection symmetry).  So L is not (globally) stable with respect to right-left equivalence since it has coinciding critical values. By general principles there should exist arbitrarily small perturbations of L which are globally stable. Plausibly, this can be achieved by moving a circle away from a symmetric position. For sufficiently small shifts, the number of stationary connecting cycles would remain the same.
\end{remark}

\noindent For $n=3$ the inradius of stationary triangle can be effectively computed.

\begin{proposition}
 The inradius $r$ of a Fermat triangle for three concentric circles with radii $a,b,c$ is the unique positive root of cubic equation
$$ 2 a b c r^3 + (a^2 b^2 +b^2 c^2 + c^2 a^2) r^2 - a^2 b^2 c^2 = 0. $$
\end{proposition}

\noindent In fact, computing the discriminant it is easy to show that this equation has only one positive real root, since the first two coefficients are positive and the free term is negative. Thus all of our aims are achieved in the case of three concentric circles. The case of four concentric circles is more subtle and exhibits new phenomena.

\section{Four Concentric Circles (4CC)}

\noindent Let $S$ be now a system of four concentric circles with pairwise different radii $r_1,r_2,r_3,r_4$.   We will consider the following types of  Fermat quadrilaterals: convex quadrilaterals, non-convex non-intersecting ones which we call \textit{spears, }eight parades which are degenerate quadrilaterals, and partially aligned quadrilaterals each of which arises from a Fermat triangle of  $S$  with one circle removed. We proceed with computing the inradius of a convex Fermat quadrilateral.

\subsection{Convex circuits 4CC}

\begin{theorem}
The inradius r of a strictly convex Fermat quadrilateral for four concentric circles with radii $a, b, c, d$ is given by:
$$
r=2 \sqrt{\frac{(Q-abc)(Q-abd)(Q-bcd)(Q-acd)}{abcd(ab+cd)(ac+bd)(ad+bc)}},
$$
where $2Q = abc + abd + acd + bcd$.
\end{theorem}

\noindent The proof uses Proposition 3 applied to properly chosen auxiliary triangles and an algebraic procedure in the spirit of elimination theory.

\begin{remark} Notice that the numerator in the above formula is not always positive so convex Fermat quadrilaterals do not exist for some collection of four radii (e.g., if  $a$ is small and $b, c, d $ are sufficiently large). As an exception we mention the partially aligned stationary circuits, which are in fact triangles. Notice also an apparent analogy with the formula for the circumradius of cyclic quadrilateral given in \cite{Var} which according to \cite{Bl} was actually already known to A.F. M\"obius.
\end{remark}

\subsection{Spears in 4CC case}

\noindent Similar results are available for spears. In the case of a spear-shaped Fermat quadrilateral with incoming second and third sides the inscribed circle should be substituted by the outscribed one, tangent to the second and third sides and tangent to the continuations of the first and fourth side. Recall that a mnemonic rule for modifying the afore mentioned formula for circumradius given in \cite{Var} was that the lengths of "negative" sides of a cyclic polygon should be taken with minus sign. Using this analogy one can produce a conjectural formula for the radius of  the tangential circle of spear by replacing some of the sides  by their negatives. We omit the details but wish to note that, like in the convex case, the expression under the square root in the conjectural formula in the spear case is not always positive, which suggests that stationary spears may not exist in some cases. 

\medskip

\noindent Knowing that convex stationary circuits and spears may not exist in some cases one can wish to study their lifespans by changing one of the given radii and creating movies. Geometrical considerations suggest that two spears may collide with a parade and create a degenerate critical point of perimeter of multiplicity three. This is illustrated by an example at the end of this section. We believe that this is a general mechanism of losing Morse property of perimeter.  The problem of characterizing collections of radii for which $L$ is not Morse remains open.

\subsection{Partially aligned for 4CC}

\noindent Partially aligned stationary circuits are solutions of the 3CC problem for any three
  out of the four circles. The remaining circle can intersect the 3CC stationary
circuit (triangle) in 2, 1 or 0 points. In the intersection points we have refraction. 
This gives respectively 2, 1 or no partially aligned solutions of the 4CC problem. 

\medskip 
\noindent
 In case of a single intersection point, the remaining circle is tangent to the triangle.
This is an intermediate case between convex and spears. It also satisfies the
 reflection rule. At this place one can expect a bifurcation where two partially 
 aligned circuits and a convex circuit transform into a spear. This is another general 
mechanism of losing the Morse property of perimeter. Obviously, both these 
 mechanisms have analogies in cases with more than four circles.

\subsection{Self-intersecting stationary circuits do not exist for 4CC}

\noindent
 We first mention that we have necessarily a Snellius circuit, otherwise we have at
 least one refraction and therefore a convex triangle. A Snellius circuit wraps
 around the socle strictly clockwise or anticlockwise over polar angles less than $\pi$. 
 Its winding number cannot therefore reach 4$\pi$, which is necessary for a tangential
self-intersecting stationary cycle.

\subsection{Parades for 4CC}

\noindent We assume now that the four radii are pairwise different.  Since all vertices are on the same straight line we assume that this line is the $x$-axis and that the fixed point on the fourth (outer) circle has positive $x$-coordinate.  Then there are eight different parades.  Let $x_i$  be the $x$-coordinate of the i-th vertex.  We have
 $x_i=\pm  r_i$   if    $i=1,2,3$ and $r_4 =r_4$. Then a direct computation shows that the Hesse matrix has the form:
$$ \left( \begin{array}{ccc}
b_1 + b_2 & - b_2             & 0 \\
- b_2            & b_2 + b_3 &  -b_3 \\
0               & - b_3             & b_3 + b_4 
\end{array}
\right)
$$
\noindent
where $b_i = \frac{x_{i-1}x_i}{|x_{i-1} - x_i|}, \;  i \; \mbox{modulo} \; 4$.
The hessian (determinant) is as follows:
$$
H(x) = \frac{x_1 x_2 x_3 x_4}{|x_1-x_2||x_2-x_3||x_3-x_4||x_4-x_1|}.
(x_1x_2|x_3-x_4|+ x_2 x_3 |x_4-x_1| + x_3 x_4 |x_1-x_2| + x_4 x_1 |x_2 - x_3|).
$$
For each parade and depending of the order of radii, this can be reduced to 
$$ 
H(x) = C(x).S(x) = C(x) . ( \frac{\epsilon_1}{r_1}   +  \frac{\epsilon_2}{r_2} +  \frac{\epsilon_3}{r_3} +  \frac{\epsilon_4}{r_4}  ),
$$
where $C(x)$  is a positive rational function and the $\epsilon_i$  take values  $-2, 0, 2$  only. These values depend on the type of the parade and on the order of the radii. We give now some illustrating examples.

\noindent  \textbf{A: }
$r_1 < r_2 < r_3 <r_4$ and $x_1 =  r_1 ; x_2 = r_2 ; x_3 =  r_3 ; x_4 =r_4$. Then
$$ S(x) = \frac{2}{r_1}+ \frac{-2}{r_4} $$
The hessian (determinant) is positive. Using Sylvester's principle it follows that this parade is a minimum.  This parade is called the \textit{shortest parade}.

\noindent  \textbf{B: }
$r_1 < r_2 < r_3 <r_4$ and $x_1 =-  r_1 ; x_2 = r_2 ; x_3 = - r_3 ; x_4 =r_4$. Then
$$ S(x) = \frac{-2}{r_1}+ \frac{-2}{r_2} +  \frac{-2}{r_3} +  \frac{-2}{r_4} $$
Note that, in this case, the hessian determinant is always negative. Using Sylvester's principle, it follows that this parade is a maximum.

\noindent  \textbf{C: }
$r_2 < r_1 < r_3 <r_4$ and $x_1 =-  r_1 ; x_2 = r_2 ; x_3 = - r_3 ; x_4 =r_4$. Then
$$ S(x) = \frac{-2}{r_1}+ \frac{2}{r_2} +  \frac{-2}{r_3} +  \frac{-2}{r_4} $$
In this case, the hessian determinant changes sign for a special set of radii,

\noindent where bifurcation can take place in situations where not all radii are pairwise different. The following example shows two spears that collide and transform into a parade.

\medskip

\noindent \textbf{Example. }We consider a symmetric situation and choose concrete numerical values of radii:
$r_1=r_3=3$ and $r_4= 4.6$
 The reader can use these values to reconstruct this example in Geogebra, Mathematica or other software. We start with  $r_2 =2.53$ and decrease  $r_2$. We follow a convex stationary circuit, which changes at $r_2= 1.7 $ (values are approximate) to a partially aligned triangular circuit which is tangent to $C_2$. At this value of $r_2$ we have a bifurcation, where two aligned circuits (which already existed for $r_2 < 1.7$  fuse with a convex solution and become a spear (due to the symmetric choice, two spears give the same image).  We decrease $r_2$ further. Around $r_2 = 1.14 $. the spear starts closing its "mouth" and about $r_2=1.13 $ the spear disappears at a zero of the Hessian. At that moment a parade which already existed for $r_2 < 1.13$ and was a saddle of index 1, fuses (pitchfork bifurcation) with the spears and  changes into a parade of index 2 which survives for smaller values of $r_2$, where there are no more bifurcations. Note that in this example will have some specific additional properties because of its symmetric nature.

\noindent 
\section{Concluding Remarks}

\noindent
 There are quite a number of natural problems suggested by the results of this note.
 For example, it is interesting to find exact conditions on four radii which
 guarantee existence of convex stationary quadrilaterals. In cases where convex 
 stationary circuits do not exist, it is interesting to find out which connecting cycle
 has the maximal perimeter. Examples suggest that this should be the longest
  parade.

\medskip

\noindent
 For $ n >4 $, there may also exist self-intersecting stationary cycles. E.g., if all radii
 are equal, the inscribed regular five-pointed star (pentagram) is a stationary
  connecting cycle. Direct computation shows that pentagram is a non-degenerate
 local maximum of perimeter. It follows that five-pointed stationary configurations 
 of local maximum type exist for non-necessarily equal but sufficiently close radii.
 By our Theorem 1 such a stationary cycle $P$ has a tangential circle inscribed in the
 convex core of P defined as the maximal convex 5-gon contained in $P$. What are the perimeters of  $P$ and its core?

\medskip

\noindent Analogous problems can be considered for other target functions on the product of 
 concentric circles. For example, one can consider stationary configurations of 
 Coulomb energy of equal charges freely sliding along the given circles. In that
case, one looks for local minima and other critical points of Coulomb potential. 
Parades are again critical for Coulomb potential. Are they generically non-degenerate critical points and what are their Morse indices? There is good evidence
 that, for Coulomb potential, there should appear all phenomena described above.

\medskip

\noindent The same problems arise for non-equal charges confined to concentric circles. By 
 a general paradigm of singularity theory, one may consider values of charges as
 parameters of the model and search for catastrophes in the sense of R.Thom. This
 seems promising since the bifurcations described in the present paper can also be t interpreted in the spirit of catastrophe theory.

\medskip

\noindent Another line of development is concerned with the same problems for disjoint
 contours with disjoint interior domains. For interior disjoint circles, it is easy to 
 obtain an analogy of Theorem 1. In this setting, an important role is played by
 mutual positions of circles. The easiest case is if there are three circles such that 
 the convex hull of any two of them does not intersect the third one. Then one can
 expect that the number of stationary circuits is eight and try to describe their geometry and calculate their Morse indices. In fact, one can easily formulate a plenty of problems in this spirit which can be successfully studied us
 the paradigms described in this paper.

\end{document}